\documentclass[preprint,12pt]{elsarticle}

\usepackage{amssymb}
\usepackage{amsmath}
\usepackage{booktabs}%
\usepackage{algorithm}%
\usepackage{algorithmicx}%
\usepackage{algpseudocode}%
\usepackage{listings}%
\usepackage{makecell}
\usepackage{multirow}%
\usepackage{enumitem}
\usepackage{xcolor}
\usepackage{longtable}
\usepackage{hyperref}
\usepackage{setspace}
\usepackage{graphicx}
\usepackage{adjustbox}
\usepackage[colorinlistoftodos]{todonotes}
\usepackage{lineno}
\usepackage{subcaption}

\usepackage{array}
\newcolumntype{L}[1]{>{\raggedright\arraybackslash}m{#1}}

\begin{document}

\begin{frontmatter}



\title{Late Train-Crew Rescheduling: A Tabu-Search-Based Approach Using Restricted DFS}


\author[a]{Liyun Yu \corref{cor1}}
\cortext[cor1]{Email: liyun.yu@liu.se}

\author[a]{Carl Henrik Häll}

\author[a]{Anders Peterson}

\author[a]{Christiane Schmidt}

\affiliation[a]{organization={Department of Science and Technology, Linköping University},
            city={Norrköping},
            country={Sweden}}

\begin{abstract}
In this paper, we reschedule the duties of train drivers one day before the operation. Because several drivers are absent (for example, because of sick leave), some trains have no driver. Thus, duties need to be rescheduled for the day of operation. We start with a feasible crew schedule for each of the remaining operating drivers, a set of unassigned tasks originally assigned to the absent drivers, and a set of standby drivers with fixed start time, end time, start depot, and end depot. Our aim is to generate a crew schedule with as few canceled tasks as possible. We present a tabu-search-based approach with restricted graph search for crew rescheduling. We also adopt a column-generation approach as a benchmark to reflect the same restrictions and objective function as our tabu-search-based approach. We use it to compare the quality of the results, computational time, and other performance indicators. Our tabu-search-based approach needs less computational time than the column-generation approach to compute an acceptable result. We further test the performance of our approach under different settings. One of the datasets used in the experiments originated from a regional passenger-train system around Stockholm, Sweden, and was provided by Mälartåg. The other dataset is generated from the open-source crew scheduling data from the ROMSOC project.
\end{abstract}


\begin{highlights}
    \item Research highlight 1: a fixed-timetable crew-rescheduling problem with driver absences and partially uncovered duties is introduced.
    \item Research highlight 2: a restricted-graph-search mechanism generates feasible alternative duties.
    \item Research highlight 3: a tabu-search-based approach is introduced to guide the neighbourhood search.
\end{highlights}


\end{frontmatter}



\section{Introduction}
Railway disruptions, caused by, for example, late changes in the railway schedule, including timetable, rolling stock, and crew schedules, have been a broadly discussed topic for many decades. While that is certainly not the preferred measure to counter disruptions, railway undertakings (RUs) often react by cancelling trains. As a consequence, passengers possibly face severe delays during travel. It is important to recover the original schedule as fast as possible after major interruptions. 

Railway undertakings aim to cancel as few trains as possible during rescheduling, as there are significant costs associated with cancellations. These costs include both financial loss and employees’ working-time loss, because employees need to spend extra time to move the canceled train back to a specific depot. Changes and complete cancellations on short notice are negative, not only for the involved services, but also for the reliability of rail transport. Furthermore, several RUs are competing in Europe. They need to gain passengers’ trust by providing a more reliable service than their competitors, for example, fewer cancellations of trains. 

Crew members are a resource, which, in case of absent crew members, must be reallocated. In this paper, we explore short-term crew rescheduling (without the need for rescheduling other resources) and try to cancel as few trains as possible. We consider the case where several train drivers take leave for the next day, and the RU needs to generate a new schedule for the remaining on-duty drivers and some standby drivers starting from various depots. We allow train cancellations as an option, but one that should be avoided. Hence, we aim to cancel as few trains as possible. 

To do so, we need to re-assign the unassigned tasks of the absent drivers to other available drivers. We do this by introducing a tabu-search-based heuristic approach and compare it to an adapted column-generation approach from Breugem et al.~\cite{breugem2022column}. Initially, the shift schedules for all on-duty drivers (both drivers with assigned tasks and standby drivers) are feasible. Because a certain number of drivers take leave at short notice, we need to re-assign the tasks that were assigned to these absent drivers to the remaining on-duty drivers---these tasks constitute a set of unassigned tasks. Thus, we are given a set of drivers who have a feasible schedule with assigned tasks in the original timetable, a set of unassigned tasks, and a set of standby drivers without any assigned tasks. 

Our goal is to develop a method suitable for handling this practical rescheduling problem. We propose a tabu-search-based approach. To obtain a feasible schedule, we allow deadheading tasks, that is, a driver using a train driven by another driver to move between stations. We limit ourselves to generating a one-day crew schedule and only consider deadheading by the trains of the considered RU, that is, no trains of other RUs can be used for deadheading. To evaluate the computational performance of our approach, we adapt a column-generation approach by \citet{breugem2022column}, such that both handle the same constraints on the resulting schedules. 

This paper is organized as follows. We present the related work in Section~\ref{Related Work}. In Section~\ref{approach}, we describe our tabu-search-based approach. Then, we run experiments and analyze the results in Section~\ref{casestudy} and conclude the paper in Section~\ref{Conclusions}. Furthermore, additional details on the adapted column-generation-approach benchmark and the performance assessment measures are provided in \ref{CG} and~\ref{eva:difficulty_weighted_assignment_rate}, respectively.

\section{Related Work} \label{Related Work}

Railway scheduling mainly involves the steps of train timetabling, rolling-stock scheduling, and crew planning. A general introduction to the field can be found in, for example, the works of \citet{cacchiani2014overview} and \citet{heil2020}. Crew planning includes two steps: crew scheduling and crew rostering. In crew scheduling, duties consisting of several tasks are formed, where a duty usually covers a short time, for example, one day \cite{abbink2018railway}. In crew rostering, duties are combined for weekly schedules that are assigned to different crew members \cite{jutte2011optimizing}. A task is characterized by several properties, including the start time, start depot, end time, and end depot. 

The common rescheduling process in railway starts with rescheduling the timetable, followed by rolling stock rescheduling, and then ends with crew rescheduling.
Crew rescheduling can target both the short term and the long term. In short-term crew rescheduling, the goal is to obtain a one-day schedule, that is, a complete working shift, including breaks; in long-term crew rescheduling, the goal is to obtain schedules of more than a day. To the best of our knowledge, most researchers consider either long-term crew rescheduling \cite{heil2020,huisman2007column}, or short-term rescheduling when the crew schedules are affected by changes in railway timetables \cite{abbink2018railway,potthoff2010,verhaegh2017,yuan2022}. In this paper, we consider short-term rescheduling that becomes necessary solely because of a change in crew availability and without underlying changes to the timetable. \citet{cacchiani2014overview} mentioned that although set-cover formulations are commonly used for both crew scheduling and rescheduling, they are not suitable for real-time rescheduling due to long computational times. 

We discuss the limitations considered among mathematical models and heuristic approaches in Subsection~\ref{limitations} and present the approaches used in both railway crew scheduling and rescheduling in Subsection~\ref{approaches}. In Subsection~\ref{contribution}, we emphasize our research goal and contributions in this paper.

\subsection{Common Restrictions}\label{limitations}

There are several common types of constraints for the schedules. The basic ones include a geographical and chronological connection between consecutive tasks and the location restrictions imposed by the start and end depots of all drivers. Other time-related constraints can be classified into two main categories. In the first category, limits on the working time, including both the total working time and uninterrupted working time, are imposed. The total working time restricts the total working time for each driver, including work, deadheading, transit, and break time \cite{boschetti2004,kokubo2017generation}. \citet{chen2013} assigned a break between two consecutive tasks and limited the maximum length of each duty. Instead of limiting the working time, the total number of tasks can be limited~\cite{froger2015}. The uninterrupted working time gives an upper time limit on consecutive work time without a break. \citet{veelenturf2016} limited the length of the working time before and after a break separately. \citet{fuentes2019} directly limited the maximum continuous driving time without any break. In the second category of time-related constraints, limits on the breaks -- including the number, the duration, and the starting-time range of breaks -- are imposed. \citet{verhaegh2017} explored different break durations and their effects.~\citet{kokubo2017generation} restricted the break duration to be within a given range. In particular, they limited the time interval in which lunch breaks may start, as did \citet{zhou2016}.

\subsection{Approaches}\label{approaches}

Due to the large size of our instances, we aim to explore the possibilities of heuristic approaches instead of exact methods. Approaches based on column generation are the most common option in crew rescheduling. \citet{huisman2007column} introduced a large-scale set covering formulation for crew rescheduling and solved it with a column-generation approach. \citet{potthoff2010} introduced another column-generation approach that can dynamically select a set of promising duties based on the approach by \citet{huisman2007column}. \citet{breugem2022column} proposed a heuristic framework combining column generation with a diving heuristic. Starting from the fractional Linear Programming (LP) solutions generated by column generation, the diving heuristic iteratively selects and fixes promising rosters, gradually constructing a feasible integer solution. Despite the numerous suggested column-generation-based approaches, \citet{heil2020} pointed out that the computational time of column generation can be very long due to problem size and slow convergence.

Additionally, \citet{heil2020} cautioned that not all meta-heuristics are a good choice to solve crew-rescheduling problems: genetic algorithms cannot maintain the feasibility of schedules in all iterations. Because in our problem, many drivers with feasible schedules of assigned tasks are present before rescheduling, we use tabu search. There are some studies of using tabu search in crew scheduling \cite{elizondo2010,chew2001,Kokubo2018}. \citet{chew2001} used tabu search for their set-partitioning formulation for crew scheduling. They introduced three move operators: exchanging one task in each duty, exchanging a sequence of tasks in each duty, and moving one task to another duty. \citet{guillermo2009} combined tabu search and integer programming to increase the diversification level and decrease the possibility of local optima. \citet{Kokubo2017train} presented a tabu-search-based approach with three types of schedule movements for generating neighborhood schedules. To make the tabu search more powerful when scheduling, they added three more categories in \citet{kokubo2017generation}. Subsequently, they changed their approach to parallel tabu search \cite{Kokubo2018}.

In most crew-rescheduling studies, disruptions caused by timetable modifications are considered, where a large proportion of the original schedules become infeasible. In contrast, we consider driver absences under a fixed timetable, where most driver schedules remain feasible, and only a subset of duties becomes uncovered. Consequently, approaches designed for timetable-disruption recovery are not directly applicable to the problem considered here, and additional functionality may yield runtimes that exceed acceptable bounds.

Among the most relevant approaches, both \citet{verhaegh2017} and \citet{yuan2022} built a depth-first iteratively-deepening search tree for searching better feasible schedules. Each node in the tree is stored with a feasible schedule and a set of remaining unscheduled tasks. A child node is generated by inserting one unscheduled task into one duty. Hence, one parent node may have a large number of child nodes. Thus, this approach is effective for small disruptions where one or a few unscheduled tasks need to be inserted. To limit the tree size, both groups of authors restrict the set of admissible connections for solving geographical conflicts: \citet{verhaegh2017} allowed at most one or two connection tasks, whereas \citet{yuan2022} defined a predefined set of connection types. Such restrictions reduce the search space by limiting allowable connections. However, they may also largely exclude feasible duties that satisfy the operational constraints but cannot be represented by the predefined connection types. Furthermore, \citet{breugem2022column} used a global search approach. Alternative rosters are searched for globally in the pricing problem, which is formulated as a resource-constrained shortest path problem in pre-built pricing graphs.

\subsection{Research Goal and Contributions}\label{contribution}

\begin{table}[b!]
\centering
\footnotesize

\renewcommand{\arraystretch}{1.2}
\setlength{\tabcolsep}{3pt}

\begin{subtable}{\columnwidth}
\centering

\begin{tabular}{@{}L{0.17\columnwidth}L{0.1\columnwidth}L{0.29\columnwidth}L{0.1\columnwidth}L{0.12\columnwidth}L{0.15\columnwidth}@{}}
\hline
Reference
& Problem type& Problem cause& Time horizon
& Feasibility scope
& Search framework
\\ \hline

\citet{breugem2022column}
& Res & Planned maintenance / construction; revised timetable and rolling-stock schedule & Multiple days
& Duty + roster
& CG-based diving
\\ \hline

\citet{Kokubo2018}
& S& Given timetable &One day
& Duty
& Parallel TS
\\ \hline

\citet{verhaegh2017}
& Res & Small disruptions: delays, cancellations, detours, or staff absences & One day
& Duty
& Iterative-deepening DFS
\\ \hline

\citet{yuan2022}
& Res & Small disruptions; revised timetable & One day
& Duty
& DFS with backtracking
\\ \hline

This paper
& Res & Driver absences; fixed timetable & One day
& Duty
& TS
\\ \hline

\end{tabular}

\caption{Overall comparison. Res stands for rescheduling problem, while S indicates scheduling problem.}
\label{tab:paper_diff_a}
\end{subtable}

\vspace{3mm}

\begin{subtable}{\columnwidth}
\centering

\begin{tabular}{@{}L{0.17\columnwidth}L{0.28\columnwidth}L{0.28\columnwidth}L{0.23\columnwidth}@{}}
\hline
Reference
& Duty-construction mechanism
& Construction restrictions
& Construction role in the search
\\ \hline

\citet{breugem2022column}
& RCSPP-based duty generation
&  Time windows; duty length; meal-break rules
& Pricing columns
\\ \hline

\citet{Kokubo2018}
& Reconstruction of one or two duties
& 6 predefined operators;\ break rules
& TS neighborhoods
\\ \hline

\citet{verhaegh2017}
& Task insertion into one duty
& Positioning trips ($\leq2$ tasks)
& DFS Tree branches
\\ \hline

\citet{yuan2022}
& Task insertion into one diagram
& 9 predefined connection types ($\leq2$ tasks)
& DFS Tree branches
\\ \hline

This paper
& DFS-based graph search for task insertion
& DFS bounds
& TS neighborhoods
\\ \hline

\end{tabular}

\caption{Duty construction comparison}
\label{tab:paper_diff_b}
\end{subtable}

\caption{Comparison with related studies}
\label{tab:paper_diff}
\vspace{-4mm}
\end{table}

In Table~\ref{tab:paper_diff}, we provide a component-level comparison between our study and the most relevant studies. We compare the problem type, problem cause, time horizon, feasibility scope, and outer search framework in Subtable~\ref{tab:paper_diff_a} and further distinguish the duty-construction mechanism, the restrictions imposed during duty construction, and the role of the constructed duties in the respective search procedures in Subtable~\ref{tab:paper_diff_b}.

For the problem setting, we consider a different problem, which is that only the train-driver schedules (but not the timetable) are impacted, and that we have feasible driver-task assignments for a (large) subset of the drivers and tasks. For us, the original timetable yields a feasible schedule for a group of drivers, and we have a set of unassigned tasks. The motivation for other crew-rescheduling problems is typically a timetable that became inapplicable because of external factors. Hence, the original crew schedule for most of the drivers is infeasible in these cases.

Our proposed method also differs substantially from the most relevant solution approaches in the literature. Existing local-search approaches generate a neighborhood through predefined insertion operators and connection types and explore this neighborhood using depth-first search (DFS) in a tree, where each leaf node represents a complete feasible schedule. In contrast, our proposed method generates feasible neighborhood candidates through restricted graph search using DFS and explores them with a tabu-search-based framework. Compared to the approaches based on global search, our proposed method uses a directed graph containing all operationally feasible task connections but bounds the alternative path exploration with restrictions of DFS.

Our main contribution can be classified into three parts. The first contribution is that we consider a problem that has not been studied so far. We investigate a railway crew-rescheduling problem caused by driver absences under a fixed timetable, where only a subset of duties becomes uncovered while most original schedules remain feasible. This is a problem that RUs face in their daily operations. The second one is that we propose a restricted-graph-search mechanism that dynamically generates feasible repair paths without predefined insertion templates or connection types. The third contribution is that we develop a tabu-search framework with tailored tabu mechanisms to explore graph-generated neighbourhoods.

\section{Tabu-Search-Based Approach}\label{approach}

In Table~\ref{tab:TB_symbol}, we summarize the symbols used in Algorithm~\ref{alg:TB}. The overall structure of our proposed tabu-search-based approach is then presented in Algorithm~\ref{alg:TB}.

\begin{table}[htb!]
\vspace{-2mm}
\centering
\begin{tabular}{@{}p{0.10\columnwidth}m{0.86\columnwidth}@{}}
\hline
Symbol & Description \\ \hline
$s_0$ & Initial feasible schedule for all remaining available drivers \\
$s^{*}$ & Historical best solution \\
$s^{\prime}$ & Current best solution \\
$T$ & Tabu list\\
$I_{\max}$ & Maximum number of iteration\\
$N(s,g)$ & Neighborhood candidates for inserting task $g$ into schedule $s$\\
$Obj(s)$ & Objective function value of schedule $s$\\ \hline
\end{tabular}
\caption{Symbols in Algorithm~\ref{alg:TB}}\label{tab:TB_symbol}
\vspace{-2mm}
\end{table}

We define the set of tasks originally assigned to the absent drivers in the original schedule as the set of unassigned tasks. The approach starts with the initial schedule, $s_0$, including the remaining drivers with feasible schedules, the set of unassigned tasks, and the set of standby drivers associated with their start and end working time and depot location.

We randomly select one unassigned task. Then we follow the structure of tabu search. We first generate neighborhood candidates by restricted graph search. To determine the best schedule among neighboring solutions, we find the one with the minimum objective value. We use a lexicographic objective that prioritizes minimizing unassigned tasks over reducing total duty length. For implementation, this is represented by a weighted sum in which the penalty for one unassigned task dominates any possible variation in duty length. This specification is mainly adopted to enable a consistent comparison with the column-generation benchmark, which uses the same objective components. The tabu-search framework, however, can accommodate alternative objective variants, including penalties for the number of drivers exceeding different overtime thresholds.

\begin{algorithm}[htb!]
\caption{Tabu-search-based approach}
\label{alg:TB}
\begin{algorithmic}[1]
\Require A feasible schedule for all remaining available drivers $s_0$ \& a set of unassigned tasks
\Ensure Historical best solution $s^{*}$
\State Initialize tabu list $T$
\State Current best solution $s^{\prime}$, Historical best solution $s^{*} \Leftarrow$ initial solution $s_0$

\For{$i = 1,\ldots,I_{\max}$}
    \State Select an unassigned task $g$
    \State Generate neighborhood candidates $N(s^{*},g)$ by restricted graph search
    \If{$N(s^{*},g)=\emptyset$}
        \State Apply diversification if required
        \If{termination condition is satisfied}
            \State \textbf{break} 
        \EndIf 
        \State \textbf{continue}
    \EndIf
    \State current best solution $s^{\prime} \Leftarrow$  the best candidate among $N(s^{*},g)$ 
    \If{Diversification criterion is satisfied \textbf{or} $Obj(s^{\prime}) < Obj(s^{*})$}
        \State $s^{*} \Leftarrow s^{\prime}$
    \EndIf
    \State Update tabu list $T$
    \If{termination condition is satisfied}
        \State \textbf{break}
    \EndIf

    \EndFor
    \State \Return $s^{*}$
    
\end{algorithmic}
\end{algorithm}

Furthermore, the tabu list contains two types. The first type stores uncovered tasks that cannot be assigned to any candidate schedule, preventing repeated unsuccessful assignment attempts. The second type stores recently accepted task-driver assignments, preventing immediate reversal of recent modifications. For the tabu tenure, we track the iteration number of each schedule in the tabu list. If the iteration number of any stored schedule exceeds the maximum limit, then we release it from the tabu list. Meanwhile, diversification is triggered after a predefined number of consecutive failure attempts. In this case, candidate solutions with equal objective values may be accepted to encourage exploration of alternative search regions. We use three termination criteria: a limit of the maximum number of iterations ($I_{\max}$), an empty set of unassigned tasks, and a limit on the maximum number of times a task can be tried to be assigned to any driver (if all unassigned tasks left have already failed a certain number of times, the algorithm terminates). We detail the process of generating neighboring solutions in Subsection~\ref{GenerateNeighbourhoodSolutions}.

\subsection{Feasibility Requirements} \label{Feasibility}

In this subsection, we introduce the feasibility requirements while generating neighboring solutions. Overall, we consider the following restrictions, which are based on what Mälartåg uses today.
\begin{enumerate}[label=(\roman*), itemsep=-2pt, topsep=1pt, leftmargin=2.5em]
  \item Temporal and geographical consistency: Avoid temporal conflicts between any two consecutive tasks of one driver's duty and ensure that consecutive tasks end and start at a matching location.
  \item Depot consistency: Each driver should start and end their duty at the driver's home depot unless this driver was relocated in the original schedule.
  \item Drivers' license: A license restriction applies to both geographical area and vehicle type.
  \item Break: Every driver who works over five hours in total must have at least one break. Each break should last at least one hour. 
  \item Rest time: Between two consecutive duties of one driver, a rest time of at least 11 hours should be ensured. Since we only replan a one-day schedule, we convert this to a restriction on the earliest start time and the latest end time for each driver.  
  \item Total working time: We allow the drivers to work limited overtime over the originally planned working time.
\end{enumerate}

\subsection{Generate Neighborhood Solutions} \label{GenerateNeighbourhoodSolutions}

\begin{table}[b!]
\centering
\begin{tabular}{@{}p{0.19\columnwidth}p{0.77\columnwidth}@{}}
\hline
Symbol & Description \\ \hline
$\hat{D}$ & Set of candidate drivers  \\
$L^{\mathrm{conflict}}$ & List of conflicting tasks\\
$g^{-}$ & Front boundary task \\
$g^{+}$ & Back boundary task \\
$P(g^{-},g)$ & Set of feasible paths connecting task $g^{-}$ to task $g$, abbreviated as $P^{-}$ \\
$P(g,g^{+})$ & Set of feasible paths connecting task $g$ to task $g^{+}$, abbreviated as $P^{+}$ \\
$s^{\mathrm{cand}}$ & Neighbouring candidate solution\\ \hline
\end{tabular}
\caption{Symbols in Algorithm~\ref{alg:neighbour_generation}}\label{tab:NS_symbol}
\end{table}

In Table~\ref{tab:NS_symbol}, we list the symbols in Algorithm~\ref{alg:neighbour_generation}. We further display the process flow of generating neighborhood candidates in Algorithm~\ref{alg:neighbour_generation}. We start with the selected task $g$ and the historical best schedule $s^*$. The candidate driver set $\hat{D}$ is obtained by filtering available drivers according to their start and end time limits. When appropriate, a subset of feasible drivers is used to reduce the search effort. 

\begin{figure}[htb!]
\centering
\includegraphics[scale=0.76]{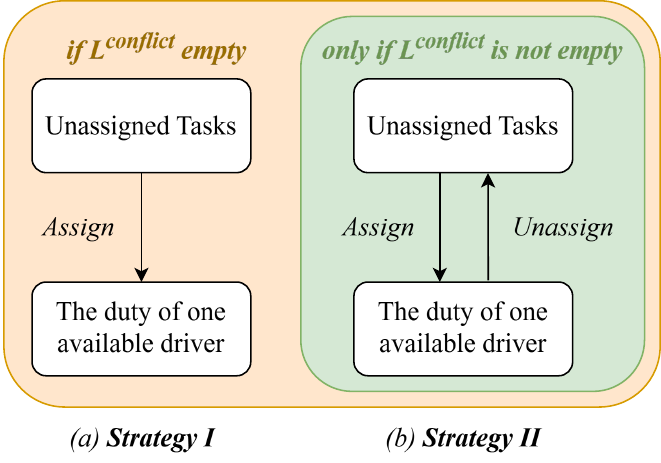}
\caption{Strategy of inserting unassigned tasks. (a): only assign unassigned tasks to one available driver; (b): assign the unassigned tasks while unassigning some assigned tasks.}
\label{fig:NeighbouringSoutionStrategies}
\vspace{-4mm}
\end{figure}

\begin{algorithm}[htb!]

\caption{Generation of neighbouring solutions}
\label{alg:neighbour_generation}
\begin{algorithmic}[1]
\Require Current solution $s^*$, selected uncovered task $g$, tabu list $T$
\Ensure Neighbourhood candidate set $N(s^*,g)$

\State $N(s^*,g) \Leftarrow \emptyset$
\For{each driver $d \in \hat{D}$}
    \State Identify $L^{\mathrm{conflict}}$ in the duty of driver $d$
    \State Initialize front boundary $g^{-}$ and back boundary $g^{+}$ based on $L^{\mathrm{conflict}}$
    
    \While{$g^{-}$ exists}
        \State Search feasible paths $P(g^{-},g)$ from front boundary $g^{-}$ to task $g$
        \If{$P(g^{-},g) \neq \emptyset$}
            \State \textbf{break}
        \EndIf
        \State Extend the front boundary $g^{-}$ one task earlier in the duty
    \EndWhile

    \If{$P(g^{-},g) = \emptyset$}
        \State \textbf{continue}
    \EndIf

    \While{$g^{+}$ exists}
        \State Search feasible paths $P(g,g^{+})$ from task $g$ to back boundary $g^{+}$
        \If{$P(g,g^{+}) \neq \emptyset$}
            \State \textbf{break}
        \EndIf
        \State Extend the back boundary $g^{+}$ one task later in the duty
    \EndWhile
    
    \For{each path in $P^{-}$ and each path in $P^{+}$}
        \State Construct $s^{\mathrm{cand}}$ by concatenating two paths
        \If{$s^{\mathrm{cand}}$ is feasible and not in tabu list $T$}
            \State $N(s^*,g) \Leftarrow N(s^*,g) \cup \{s^{\mathrm{cand}}\}$
        \EndIf
    \EndFor
\EndFor

\State \Return $N(s^*,g)$
\end{algorithmic}
\end{algorithm}

Inserting $g$ into any duty may create temporal conflicts. We define the set of tasks, which is a subset of the selected duty in $s^*$, that become infeasible while inserting $g$ as $L^{\mathrm{conflict}}$. Furthermore, we employ two strategies, which are shown in Figure~\ref{fig:NeighbouringSoutionStrategies}. For Strategy I, we only assign tasks to the drivers without deleting any tasks from the drivers' schedules (that is, without unassigning any tasks), see Figure~\ref{fig:NeighbouringSoutionStrategies}(a). For Strategy II, we try to assign tasks to a driver and delete some tasks from this driver's schedule, see Figure~\ref{fig:NeighbouringSoutionStrategies}(b). When $L^{\mathrm{conflict}}=\emptyset$, both strategies are applicable. The front boundary task $g^{-}$ and back boundary task $g^{+}$ denote the tasks immediately before and after the conflict task list in the duty list, respectively.
while only Strategy II can be applied when $L^{\mathrm{conflict}}\neq\emptyset$. In this case, $g^{-}$ is the last task ending before the start time of $g$, and $g^{+}$ is the first task starting after the end time of $g$.

The aim of the graph search is to identify feasible replacement paths that reconnect the duty between the two boundaries, passing through the selected task $g$. With restricted graph search, we search for the feasible paths on a pre-defined graph. The graph is constructed in a preprocessing step. Each node represents a task, while directed arcs are created only between pairs of tasks that satisfy both temporal and geographical consistency. The path search is executed by a depth-first search (DFS) with restrictions on the search depth and performed in two stages. First, feasible paths $P(g^{-},g)$ from the front boundary task $g^{-}$ to task $g$ are searched. If no feasible path is found, the front boundary is extended one task earlier in the duty, and the search is repeated. If the front boundary has been extended to the beginning of the duty, the search is then initialized from the driver's start depot and explores feasible paths to task $g$. Only tasks starting from the depot and also satisfying the driver's earliest start-time constraint and ending before the start of $g$ are considered during the search. This process continues until either at least one feasible path is found or the front boundary exceeds the targeted duty. Similarly, feasible paths from task $g$ to the back boundary task $g^{+}$ are searched. If no feasible path is found, the back boundary is extended one task later in the duty, and the search is repeated. Once the end of the duty is reached, the search is then performed from $g$ to the driver's end depot. When feasible paths exist for both search directions, they are concatenated to construct candidate replacement duties.

\subsection{Neighborhood Example}\label{Neighborhood_Example}

\begin{figure}[h!]
\vspace{-4mm}
\centering
\includegraphics[scale=0.71]{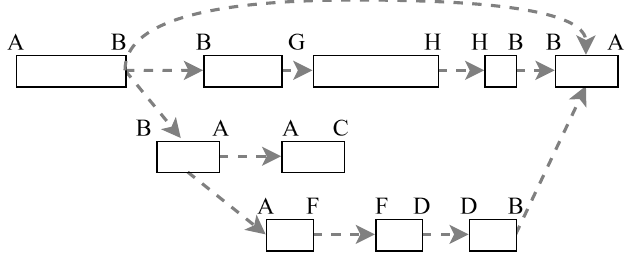}
\caption{Example of a search graph}
    \label{fig:method_example_GRAPH}
\vspace{-4mm}
\end{figure}

In Figure~\ref{fig:method_example_GRAPH}, we illustrate an example of the search graph used in the proposed approach.
Each node represents a task and is depicted as a block in Figure~\ref{fig:method_example_GRAPH}, with time shown on the horizontal axis and the block length indicating the task duration. Each dashed directed arc corresponds to a temporally and geographically feasible connection. Feasible paths are subsequently generated by performing a depth-first search within the restricted search region defined by the start and end nodes.

\begin{figure}[h!]
\centering
\includegraphics[scale=0.61]{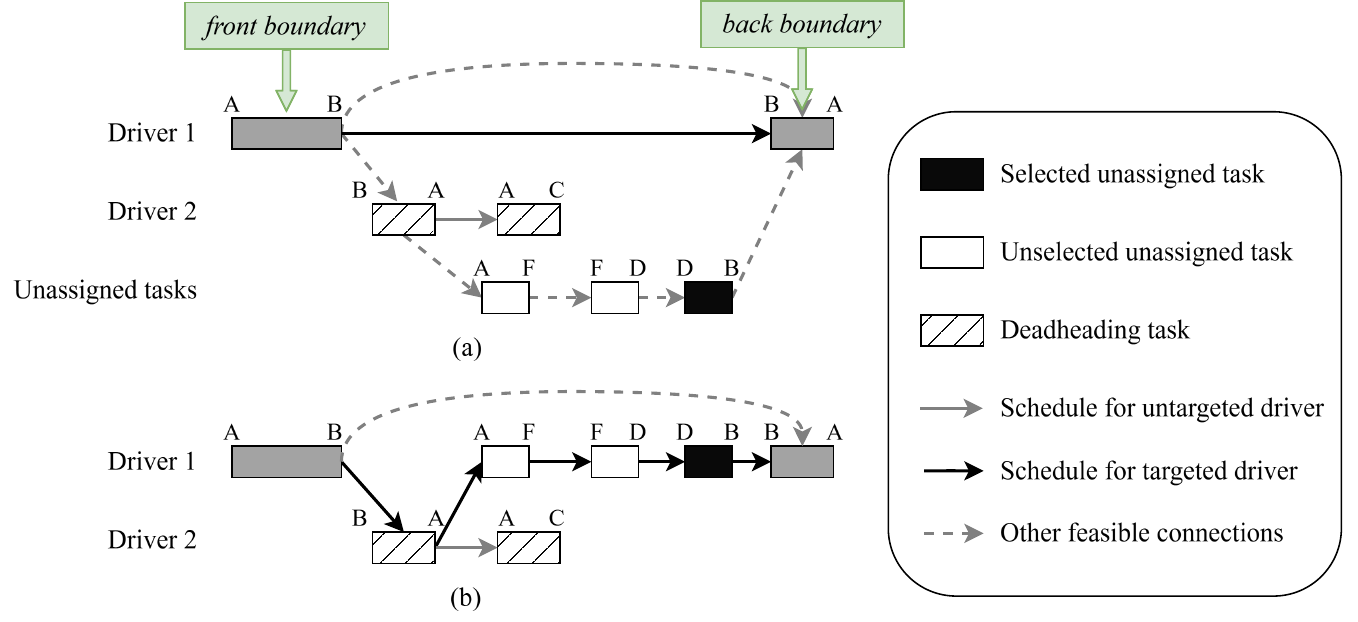}
\caption{Example of Strategy I. (a): initial schedule; (b): new schedule with Strategy I}
    \label{fig:4}
\end{figure}

\begin{figure}[h!]
\centering
\includegraphics[scale=0.61]{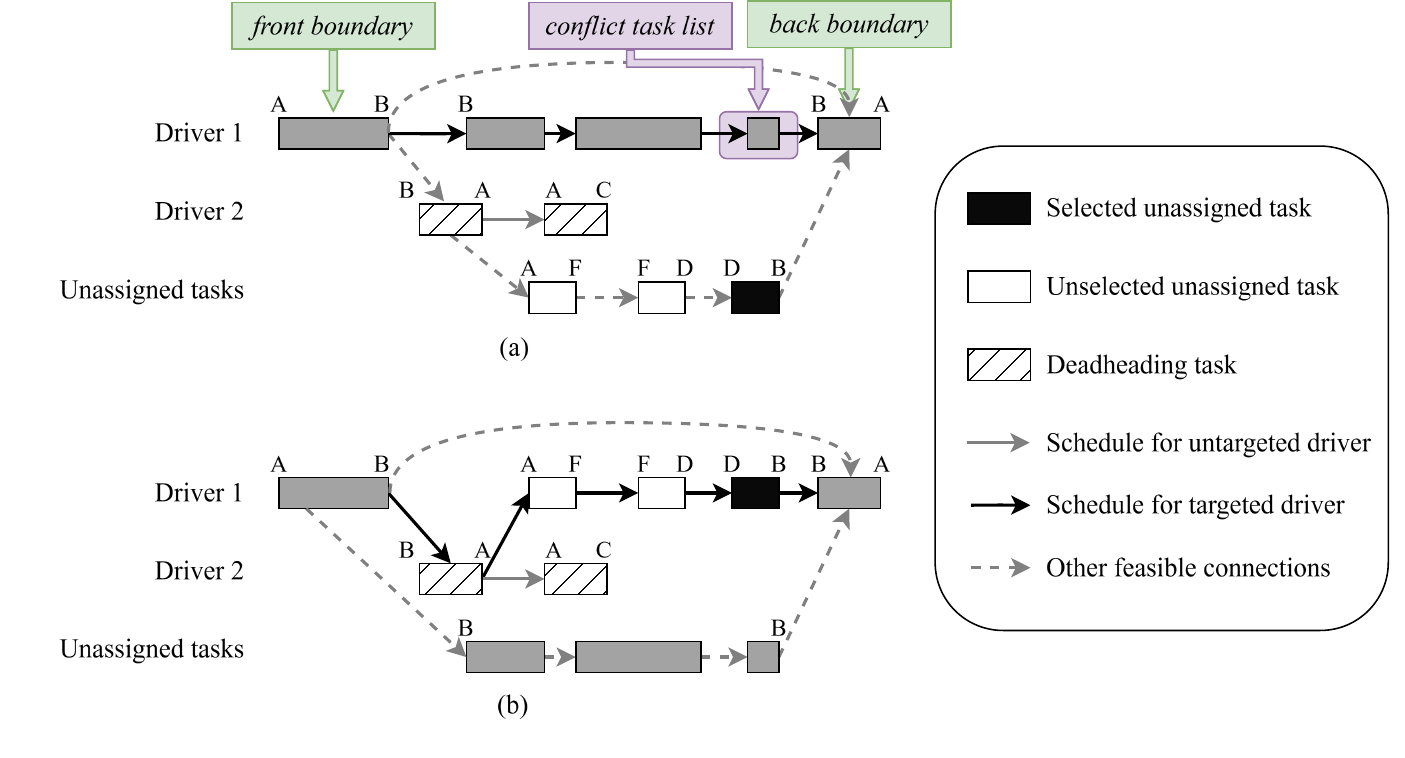}
\caption{Example of Strategy II. (a): initial schedule; (b): new schedule with Strategy II}
    \label{fig:5}
\vspace{-6mm}
\end{figure}

In Figures~\ref{fig:4} and~\ref{fig:5}, we present examples of Strategy I and Strategy II, respectively. Figure~\ref{fig:4}(a) shows the initial schedule of Driver 1, Driver 2, and the pool of unassigned tasks. The goal is to assign the selected unassigned task (the black block indicating a task with the temporal length of the block's width that runs from Station D to Station B) to Driver 1. Because Driver 1 does not have a task ending at Station D, we insert the combination of one deadheading task from the schedule of Driver 2 and two other unassigned tasks before the selected unassigned task. The resulting schedule is shown in Figure~\ref{fig:4}(b). In Figure~\ref{fig:5}(a), we try to assign the same selected unassigned task to Driver 1. However, Driver 1 has a different schedule than in Figure~\ref{fig:4}. Since there are temporal conflicts between the selected task and the initial schedule of Driver 1, we test unassigning three tasks from the schedule of Driver 1, which yields the schedule shown in Figure~\ref{fig:5}(b).

\section{Numerical Experiments}\label{casestudy}

\begin{figure}[htb!]
\centering
\includegraphics[scale=0.58]{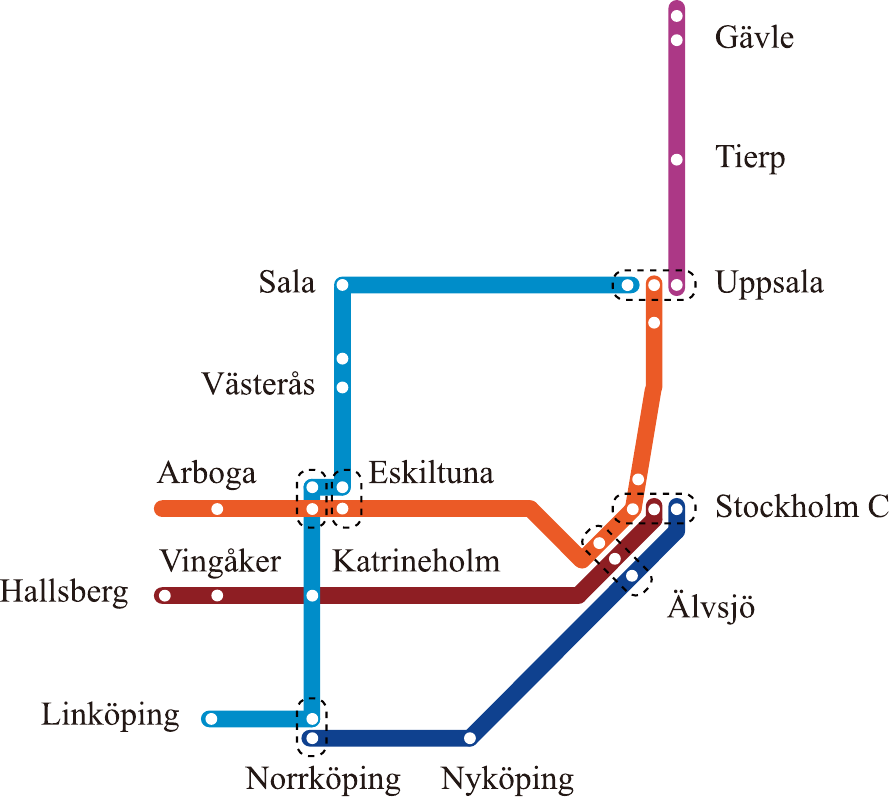}
\caption{A map of Mälartåg's service network based on their complete network map \cite{MTR}}
    \label{fig:map}
\vspace{-2mm}
\end{figure}

The numerical experiments are primarily based on a real-world dataset from a regional passenger-train system around Stockholm, Sweden, provided by Mälartåg. To complement this case study, we also use datasets derived from the open-source crew-scheduling data of the ROMSOC project~\citep{staszek_2021_5171817}. As the ROMSOC data do not include realized schedules, we first generate feasible initial schedules using the proposed tabu-search-based method and use them as input to the rescheduling experiments. For the datasets from Mälartåg, we present the RU's network in Figure~\ref{fig:map}: the traffic consists of five lines, each with an end-to-end travel time of one to three hours. The original Mälartåg network consists of 45 stations. Several of these stations are not appropriate to start or end a task. Hence, in our experiments, we only include 20 stations, and we allow drivers to swap vehicles in all of these stations marked in Figure~\ref{fig:map}. For the datasets generated based on the ROMSOC-project data, there are 42 stations considered in the railway network.

We run all the experiments on a computer with 12th Gen Intel(R) Core(TM) i5-12400F 2.50 GHz and 16.0GB RAM. We use Gurobi 11.0 to solve the restricted master problems in the column-generation approach described in Appendix~\ref{CG}. As an exact benchmark, we also solve the corresponding full master problem directly with Gurobi 11.0, in which all feasible duties are enumerated in advance and included explicitly. Our implementation of the tabu search algorithm was carried out in Python.

In Subsection~\ref{Instances}, we present details of the instances. In Subsection~\ref{subsec:result_com}, we compare the results computed by our tabu-search-based approach with the column-generation approach and the exact solver, i.e., Gurobi. Furthermore, we evaluate the graph search used in our tabu-search approach in Subsection~\ref{subsec:Graph_Search_Eva}. Finally, in Subsection~\ref{analysis}, we analyze the performance of our tabu-search-based approach under different parameter settings.

\subsection{Instances}\label{Instances}

\begin{table}[b!]
\vspace{-4mm}
\centering 
\begin{tabular}{@{}c|l|cc|cc@{}}
\hline
 \multirow{3}{*}{\makecell[c]{Instance\\source}} & \multirow{3}{*}{\makecell[l]{Instance\\name}}  & \multicolumn{2}{c|}{Number of drivers}  & \multicolumn{2}{c}{Graph size}\\ \cline{3-6}
 &  &  \multirow{2}{*}{\makecell[c]{\# operating\\drivers}}  & \multirow{2}{*}{\makecell[c]{\# standby\\drivers}}   & \multirow{2}{*}{\# nodes} & \multirow{2}{*}{\# arcs}\\ 
 &&&&&\\ \hline

 \multirow{9}{*}{\makecell[c]{Mälartåg}} & 18\_small  & 38 & 4 & 212&2191\\
                            & 18\_med     & 76 & 11 &360&6458\\ 
                            & 18\_large      & 99 & 16  &546&15043\\ \cline{2-6}
 &19\_small    & 39 & 4 &180&1768\\
                            & 19\_med     & 80 & 10 &360&7011\\ 
                            & 19\_large      & 99 & 16 &439&10987\\ \cline{2-6}
  &20\_small    & 42 & 7 &214&2123\\
                            & 20\_med   & 60 & 11 &364&6593\\ 
                            & 20\_large     & 100 & 16 &552&15123\\ \hline
\multirow{10}{*}{\makecell[c]{ROMSOC}} & 2W\_1\_1   & \multicolumn{2}{c|}{217} & 163 & 4849\\
                            & 2W\_1\_2    &\multicolumn{2}{c|}{217}  &192&722\\ 
                            & 2W\_1\_3   &\multicolumn{2}{c|}{217}   &188&768\\ 
                            & 2W\_1\_4   & \multicolumn{2}{c|}{217}   &194&729\\ 
                            & 2W\_1\_5   &\multicolumn{2}{c|}{217}   &182&702\\ 
                            & 2W\_1\_6   &\multicolumn{2}{c|}{217}  &183&804\\ 
                            & 2W\_1\_7   & \multicolumn{2}{c|}{217}   &137&463\\ 
                            \cline{2-6} 
                            & 3W\_2\_1   & \multicolumn{2}{c|}{217}  &248&1429\\
                            & 3W\_2\_2    &  \multicolumn{2}{c|}{217} & 290&1703\\ 
                            & 3W\_2\_3   & \multicolumn{2}{c|}{217}   &283&1756\\  \hline
\end{tabular}
\caption{Instances and graph sizes}\label{tab:InstanceSize_graphsize}
\vspace{-7mm}
\end{table}

The instances from Mälartåg stem from three working days (April 18\textsuperscript{th}, 19\textsuperscript{th}, and 20\textsuperscript{th} in 2023). Except for the instances with the one-day schedule for all train drivers for these days (denoted as large instances), we use subsets of all train-driver shifts for these days to create medium and small instances. We present the resulting nine instances from Mälartåg in Table~\ref{tab:InstanceSize_graphsize}. Furthermore, all ROMSOC-derived instances are constructed with a one-day planning horizon. In Table~\ref{tab:InstanceSize_graphsize}, we display all experimental instances in terms of available drivers and the size of the search graph. For the Mälartåg instances, driver resources are reported separately for operating and standby drivers, whereas for the ROMSOC-derived instances, we report the total number of available drivers. The size of the task graph is described by the number of nodes and directed arcs.

\subsection{Result Comparison}\label{subsec:result_com}
We compare the proposed tabu-search-based method with exact and column-generation benchmarks where applicable. The adapted formulation of the column-generation approach is presented in \ref{CG}. The objective is lexicographic: the number of unassigned tasks is minimized first, followed by total duty length. It is implemented as a weighted sum, assigning a penalty coefficient of 1 to each unassigned task and $10^{-6}$ to duty length measured in minutes. This penalty exceeds the maximum possible variation in duty length, thereby preserving the lexicographic priority. Due to the randomized task-selection step, our tabu-search-based approach is run 100 times for each instance, and aggregated results are reported.

\subsubsection{Comparison on Reduced-Scale Instances}\label{subsec:result-small_medium}
In this subsection, we evaluate solution quality on reduced-scale instances for which the full master problem can be solved to optimality by Gurobi. Unlike the column-generation approach, all feasible duties are enumerated in advance and included explicitly in the model. The tabu-search-based method is compared with both the exact solution and the column-generation approach.
\begin{table}[htb!]
\centering 
{\begin{tabular}{c|c|ccc}
\hline
\makecell[c]{Instance\\name} & Approach
& \makecell[c]{\# Unassigned\\tasks}
& \makecell[c]{Assignment\\rate}
& Time \\
\hline
\multirow{4}{*}{18\_small} & Exact solver & 4 & 69.2\% & 6.47 s \\
& CG & 4 & 69.2\% & 14.0 s \\
& TS [\textit{max}] & 4 & 69.2\% & 0.2 s \\
& TS [\textit{med}] & 4 & 69.2\% & 0.2 s \\
\hline
\multirow{4}{*}{19\_small} & Exact solver & 5 & 58.3\% & 2.57 s \\
& CG & 7 & 41.7\% & 11.7 s \\
& TS [\textit{max}] & 6 & 50.0\% & 0.2 s \\
& TS [\textit{med}] & 6 & 50.0\% & 0.2 s \\
\hline
\multirow{4}{*}{20\_small} & Exact solver & 5 & 50.0\% & 7.95 s \\
& CG & 6 & 40.0\% & 13.7 s \\
& TS [\textit{max}] & 5 & 50.0\% & 0.2 s \\
& TS [\textit{med}] & 5 & 50.0\% & 0.2 s \\
\hline
\end{tabular}}
\caption{Results for small-sized instances}
\label{tab:small}
\end{table}

\begin{table}[htb!]
\centering 
{\begin{tabular}{c|c|ccc}
\hline
\makecell[c]{Instance\\name} & Approach
& \makecell[c]{\# Unassigned\\tasks}
& \makecell[c]{Assignment\\rate}
& Time \\
\hline
\multirow{4}{*}{18\_med} & Exact solver & 7 & 69.6\% & 386 s \\
& CG & 10 & 56.5\% & 249.8 s \\
& TS [\textit{max}] & 8 & 65.2\% & 1.3 s \\
& TS [\textit{med}] & 8 & 65.2\% & 1.3 s \\
\hline
\multirow{4}{*}{19\_med} & Exact solver & 5 & 75.0\% & 439 s \\
& CG & 5 & 75.0\% & 333.4 s \\
& TS [\textit{max}] & 5 & 75.0\% & 1.3 s \\
& TS [\textit{med}] & 5 & 75.0\% & 1.3 s \\
\hline
\multirow{4}{*}{20\_med} & Exact solver & 15 & 76.6\% & 433 s \\
& CG & 15 & 76.6\% & 199.1 s \\
& TS [\textit{max}] & 15 & 76.6\% & 13.2 s \\
& TS [\textit{med}] & 15 & 76.6\% & 13.2 s \\
\hline
\end{tabular}}
\caption{Results for medium-sized instances}
\label{tab:medium}
\end{table}

In Table~\ref{tab:small}, we display the results of all small-sized instances, while in Table~\ref{tab:medium}, we show the results of all medium-sized instances. We present the number of unassigned tasks (that is, the number of tasks that remain unassigned at the end of the run), the assignment rate (that is, the percentage of initially unassigned tasks that we could assign), and the computational time in seconds. CG refers to the results of the adapted column-generation approach, and TS refers to the results of our tabu-search-based approach. To ensure statistical reliability, we perform 100 repetitions with our tabu-search-based approach for each small- and medium-sized problem instance, utilizing the same set of absent drivers as in the benchmark tests. We display the maximum and median assignment rates among the 100 runs where we break ties by choosing the instance with the shortest computational time. In both Tables~\ref{tab:small} and~\ref{tab:medium}, the tabu-search-based method achieves competitive solution quality while requiring less computational time than the column-generation approach. In the small-sized instances, the runtime of the tabu search is around 0.2 seconds, compared with more than 10 seconds for column generation. In the medium instances, the tabu search requires 1.3--13.2 seconds, whereas column generation requires 3--6 minutes. The tabu search also matches or closely approaches the optimal number of unassigned tasks in most cases, showing that its computational efficiency is achieved without a substantial loss of solution quality.

\subsubsection{Comparison on Full-Scale Instances}\label{subsec:result-large}
In this subsection, we compare our tabu-search-based method with the column-generation approach on the full-scale instances with a one-day horizon, where direct solution of the complete model by an exact solver is computationally intractable. For each Mälartåg and ROMSOC-derived instance, we randomly select 11 and 22 absent drivers, respectively, and solve each instance 100 times using the same set of absent drivers.

\begin{table}[htb!]
\centering 
{\begin{tabular}{c|c|ccc}
\hline
\makecell[c]{Instance\\name} & Approach
& \makecell[c]{\# Unassigned\\tasks}
& \makecell[c]{Assignment\\rate}
& Time \\
\hline
\multirow{3}{*}{18\_large}
& CG & 9 & 76.9\% & 1073 s \\
& TS [\textit{max}] & 8 & 79.5\% & 4.4 s \\
& TS [\textit{med}] & 10 & 74.4\% & 13.6 s \\
\hline
\multirow{3}{*}{19\_large}
& CG & 7 & 73.1\% & 933 s \\
& TS [\textit{max}] & 5 & 80.8\% & 1.6 s \\
& TS [\textit{med}] & 7 & 73.1\% & 2.3 s \\
\hline
\multirow{3}{*}{20\_large}
& CG & 11 & 71.8\% & 1528 s \\
& TS [\textit{max}] & 8 & 79.5\% & 4.9 s \\
& TS [\textit{med}] & 11 & 71.8\% & 12.9 s \\
\hline
\end{tabular}}
\caption{Results for large-sized instances from Mälartåg}
\label{tab:large_MTR}
\end{table}

\begin{table}[htb!]
\centering
{\begin{tabular}{c|c|ccc}
\hline
\makecell[c]{Instance\\name} & Approach
& \makecell[c]{\# Unassigned\\tasks}
& \makecell[c]{Assignment\\rate}
& Time \\
\hline
\multirow{3}{*}{2W\_1\_1}
& CG & 23 & 30.3\% &  64.5 s \\
& TS [\textit{max}] & 6 & 81.8\% & 2.6 s \\
& TS [\textit{med}] & 6 & 81.8\% & 2.6 s \\
\hline
\multirow{3}{*}{2W\_1\_2}
& CG & 24 & 20.0\% & 83.4 s \\
& TS [\textit{max}] & 9 & 70.0\% & 2.6 s \\
& TS [\textit{med}] & 9 & 70.0\% & 2.6 s \\ 
\hline
\multirow{3}{*}{2W\_1\_3}
& CG & 21 & 38.2\% & 95.0 s \\
& TS [\textit{max}] & 14 & 58.8\% & 3.7 s \\
& TS [\textit{med}] & 14 & 58.8\% & 3.7 s \\
\hline
\multirow{3}{*}{2W\_1\_4}
& CG & 28 & 17.6\% & 113.3 s \\
& TS [\textit{max}] & 5 & 85.3\% & 2.6 s \\
& TS [\textit{med}] & 5 & 85.3\% & 2.6 s \\ 
\hline
\multirow{3}{*}{2W\_1\_5}
& CG & 24 & 25.0\% & 90.6 s \\
& TS [\textit{max}] & 8 & 75.0\% & 2.7 s \\
& TS [\textit{med}] & 8 & 75.0\% &  2.7 s \\
\hline
\multirow{3}{*}{2W\_1\_6}
& CG & 27 & 27.0\% & 101.7 s \\
& TS [\textit{max}] & 10 & 73.3\% & 3.2 s \\
& TS [\textit{med}] & 10 & 73.3\% & 3.2 s \\
\hline
\multirow{3}{*}{2W\_1\_7}
& CG & 28 & 6.7\% &  59.5 s \\
& TS [\textit{max}] & 0 & 100.0\% & 1.3 s \\
& TS [\textit{med}] & 0 & 100.0\% &  1.3 s \\
\hline
\multirow{3}{*}{3W\_2\_1}
& CG & 23 & 30.3\% & 169.1 s \\
& TS [\textit{max}] & 9 & 72.7\% & 3.5 s \\
& TS [\textit{med}] & 9 & 72.7\% &  3.5 s \\
\hline
\multirow{3}{*}{3W\_2\_2}
& CG & 31 & 0.0\% & 361.0 s \\
& TS [\textit{max}] & 14 & 54.8\% & 4.5 s \\
& TS [\textit{med}] & 14 & 54.8\% & 4.5 s \\ 
\hline
\multirow{3}{*}{3W\_2\_3}
& CG & 31 & 11.4\% & 355.6 s \\
& TS [\textit{max}] & 12 & 65.7\% & 4.2 s \\
& TS [\textit{med}] & 12 & 65.7\% & 4.2 s \\ 
\hline
\end{tabular}}
\caption{Results for large-sized instances from ROMSOC}
\label{tab:large_ROMSOC}
\end{table}

In Table~\ref{tab:large_MTR}, we give the summary of test-run results for each large instance from Mälartåg. As in Subsection~\ref{subsec:result-small_medium}, we present the number of unassigned tasks, the assignment rate, and the computational time in seconds. For the tabu-search-based approach, the maximum assignment rate is in the interval $[79.5\%, 80.8\%]$, and the median rate is in the interval $[71.8\%, 74.4\%]$. Overall, the computational time is less than 20 seconds. Compared with the column-generation approach, the tabu-search-based method yields a similar median value of the number of unassigned tasks, while substantially reducing the computational time. The column-generation approach requires approximately 15--25 minutes to solve these instances.

Furthermore, we display the result comparison for ROMSOC instances in Table~\ref{tab:large_ROMSOC}. In general, our tabu-search-based method obtains much higher assignment rates than the column-generation approach with a diving heuristic in all tested instances, while requiring substantially shorter computational times. The median runtimes remain between 1.3 and 4.5 seconds, whereas the runtimes of CG range from 1 to 6 minutes. The low assignment rates of CG on these instances are mainly due to the diving heuristic used for ensuring integer solutions in the used column-generation approach. An unfavorable early fixing decision can restrict the remaining recovery options. This effect is more pronounced since the schedules generated from ROMSOC are relatively sparse. Thus, the fractional LP solution may provide limited guidance for sequential fixing.

\begin{figure}[b!]
\centering
\includegraphics[width=0.87\textwidth]{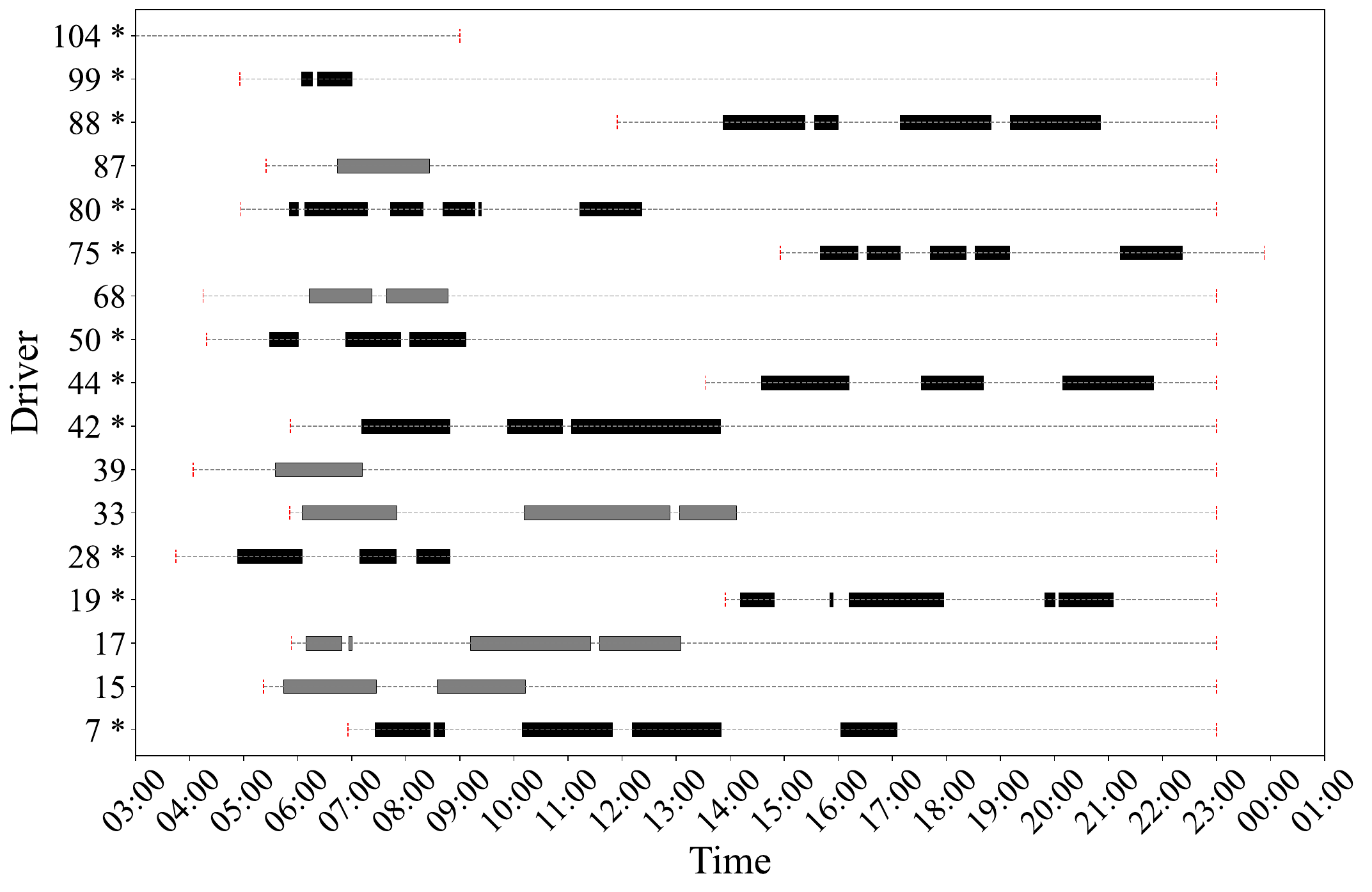}
\caption{The original schedule of the large instance on April 20th. Gray blocks indicate tasks of drivers with changed duties (in the new schedule), black blocks indicate tasks assigned to now absent drivers (which need to be reassigned in the new schedule). The start and end times of the drivers' possible shifts are shown by red markers. The absent drivers are marked with asterisks.}
\label{fig:OriginalSchedule}
\end{figure}

\begin{figure}[htb!]
\centering
\includegraphics[width=0.87\textwidth]{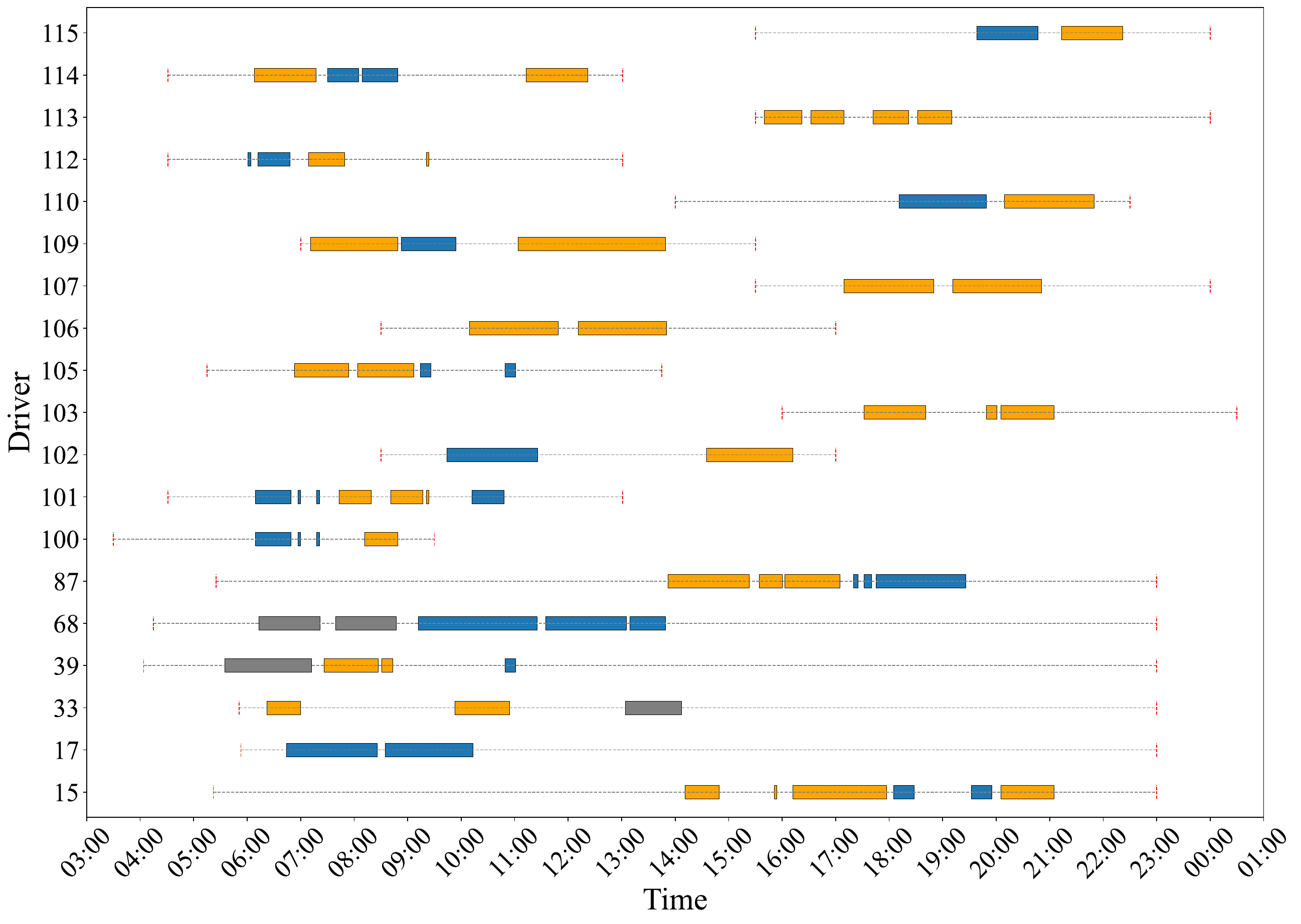}
\caption{A new schedule of the large instance on April 20th obtained with our tabu-search-based approach. Orange blocks indicate tasks that were originally assigned to absent drivers, gray blocks indicate tasks that are assigned to the same driver as in the original schedule, and blue blocks indicate tasks that were assigned to operating drivers in the original schedule but were reassigned in the new schedule. Again, the start and end times of the drivers' possible shifts are shown by red markers. All drivers with numbers $\geq 99$ are standby drivers.}
\label{fig:NewSchedule}
\end{figure}

We further present the details of one result for the large instance on April 20th. The original schedule contains 100 operating drivers and 16 standby drivers. We have 11 absent drivers who were originally operating drivers. In this run, the randomly selected 11 drivers were assigned to 39 tasks in total. Figure~\ref{fig:OriginalSchedule} shows the original schedule for all the absent drivers and the drivers with changed duties in the new schedule, but not the schedules of the standby drivers and those of drivers without a changed duty in the new schedule. Figure~\ref{fig:NewSchedule} displays all changed duties: those of on-duty and standby drivers. The 100 operating drivers in the original schedule are assigned numbers 0 to 99; drivers with a number $\geq 100$ are standby drivers. Out of the 39 tasks that we aim to reassign, in this run, we assign 32 tasks, which yields an assignment rate of $\sim$82.1 percent. Seven tasks remain unassigned in this run.

\subsection{Graph Search Evaluation}\label{subsec:Graph_Search_Eva}

In this subsection, we evaluate the restricted DFS for our graph search. For each graph, we randomly sample 1000 task pairs and search for feasible connection paths under different search limits. In Table~\ref{tab:graph_search_evaluation}, we display both the search limits and corresponding output. The search limits are defined by the maximum path length, the maximum number of paths returned for each pair, and the maximum number of visited nodes. The column
\textit{feasible pair} denotes the number of pairs for which at least one feasible path is found. The column \textit{short-path pairs} stands for the number of pairs with at least one path containing at most two nodes, i.e., tasks. The two-task threshold provides a comparison with connection paths involving at most two tasks, as considered in related studies \cite{verhaegh2017,yuan2022}. \textit{Extended-only rate} is the proportion of feasible pairs for which all identified feasible paths contain more than two nodes.

\begin{table}[htb!]
\centering
{\begin{tabular}{@{}c|ccc|ccc@{}}
\hline 
\multirow{3}{*}{\makecell[c]{Instance\\name}}
& \multicolumn{3}{c|}{Search limits}
& \multicolumn{3}{c}{Output}\\
\cline{2-7}
& \makecell[c]{path\\length}
& \# paths
& \makecell[c]{visited\\nodes}
& \makecell[c]{feasible\\pairs}
& \makecell[c]{short-path\\pairs}
& \makecell[c]{extended-only\\rate}\\ \hline
18\_large &8&3&41 & 140 & 26 & 81.4\% \\ \hline
19\_large &8&3&35 & 111 & 19 & 82.9\% \\ \hline
20\_large &8&3&41 & 122 & 15 & 87.7\% \\ \hline
18\_large &10&6&82 & 193 & 28 & 85.5\% \\ \hline
19\_large &10&6&71 & 184 & 22 & 88.0\% \\ \hline
20\_large &10&6&83 & 185 & 18 & 90.3\% \\ \hline
18\_large &12&9&206 & 254 & 33 & 87.0\% \\ \hline
19\_large &12&9&179 & 229 & 28 & 87.8\% \\ \hline
20\_large &12&9&207 & 223 & 20 & 91.0\% \\ \hline

\end{tabular}}
\caption{Evaluation of graph search with restricted DFS}
\label{tab:graph_search_evaluation}
\end{table}

Overall, the results show that a large proportion of feasible connections are extended-only. This indicates that the restricted DFS procedure identifies many feasible paths beyond those containing at most two task nodes.

\begin{table}[htb!]
\centering 
{\begin{tabular}{c|c|c}
\hline 
Instance name
&Maximum path-length limits
& Assignment rates\\ \hline
18\_large[max] & 2 & 64.1\% \\
18\_large[max] & 8 & 79.5\% \\ 
18\_large[med] & 2 & 64.1\% \\
18\_large[med] & 8 & 74.4\% \\
\hline
19\_large[max] & 2 & 61.5\% \\ 
19\_large[max] & 8 & 80.8\% \\
19\_large[med] & 2 & 57.7\% \\
19\_large[med] & 8 & 73.1\%\\ 
\hline
20\_large[max] & 2 & 71.8\% \\ 
20\_large[max] & 8 & 79.5\% \\
20\_large[med] & 2 & 69.2\% \\
20\_large[med] & 8 & 71.8\% \\ \hline
\end{tabular}}
\caption{Comparison of TB results under different limits}
\label{tab:graph_search_evaluation_2_limit}
\end{table}

In Table~\ref{tab:graph_search_evaluation_2_limit}, we further investigate the effect of altering the maximum path-length restriction based on the results from our tabu-search-based approach. The maximum path-length limit is set to eight throughout the computational study, except in this table. The results show that increasing the limit from two to eight improves the assignment rate in every reported case, with gains ranging from 2.6\% to 19.3\%. With a maximum path-length limit of two, the computational time ranges from 2 to 5~s. Increasing the limit to eight leads this range to approximately 2--20~s. The results therefore indicate that paths containing more than two connecting tasks can provide additional feasible reassignment opportunities.

\subsection{Analysis}\label{analysis}
In this subsection, we evaluate the proposed tabu-search-based method under two different settings. We first adjust the number of absent drivers. Then, we examine the effect of enforcing drivers' start and end time limits. The purpose is to assess how these settings affect the solution. In addition to the assignment rate, we also report a difficulty-weighted assignment rate (DWAR) to assess the difficulty of assigned and remaining unassigned tasks. The assignment rate treats all tasks equally, whereas the DWAR assigns task-specific weights to reflect the difficulty of assigning a task. The weight captures two factors: the structural connectivity of a task in the rescheduling graph, measured by the corresponding node's in- and out- degree, and its temporal flexibility within the planning horizon, measured by its distance to the start or end of the planning horizon, since tasks located closer to either boundary typically have fewer available drivers and thus reduced flexibility. These two components are both important for capturing the difficulty of assigning a task. To further quantify this difference, we introduce $\Delta$, defined as the difference between DWAR and the assignment rate. This is useful because the assignment rate alone does not reflect differences in task difficulty, and two solutions with the same assignment rate may still end up with unassigned tasks that differ in how difficult they are to recover. A negative $\Delta$ indicates that the remaining tasks are more difficult to assign; a positive $\Delta$ indicates that they are easier than the already assigned tasks. Detailed definitions are provided in \ref{eva:difficulty_weighted_assignment_rate}.

\subsubsection{Varying Numbers of Absent Drivers}\label{subsubsec:Tabu_Search_Eva_Nr_Driver}

\begin{figure}[htb!]
\vspace{-4mm}
\centering
\includegraphics[scale=0.3]{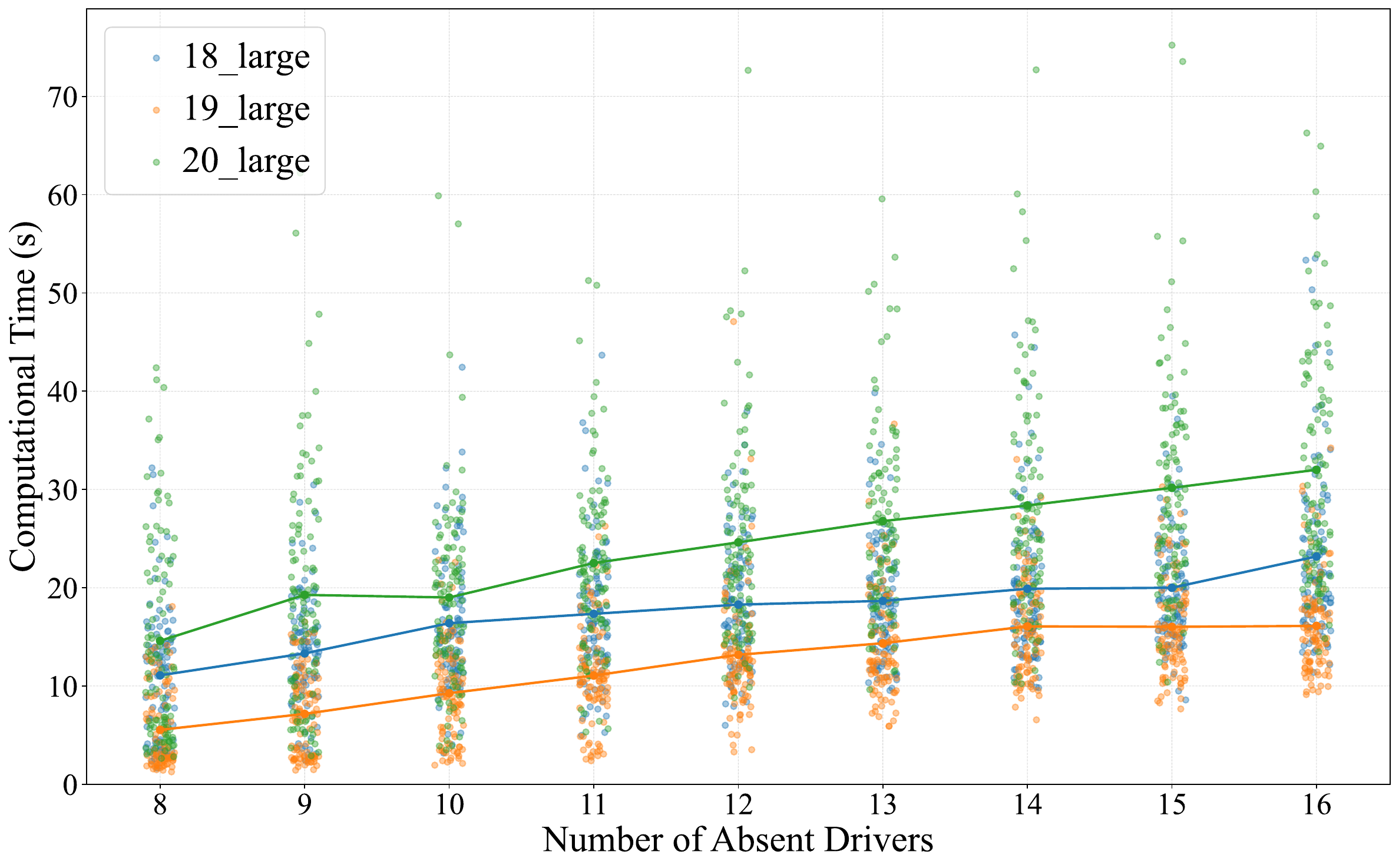}
\caption{Computational time of the large instances under different numbers of absent drivers. Points represent individual runs. Solid lines connect the median values computed for each number of absent drivers.}
\label{fig:ComputationTime_NrAbsentDrivers}
\vspace{-4mm}
\end{figure}

We vary the number of absent drivers in the large instances. In our real-world scenario, typically 8 to 11 drivers are absent per day. Here, we vary the number of absent drivers between 8 and 16. For each number of absent drivers in each instance, we run the heuristic approach 100 times. For each run, we randomly select the absent drivers. Figure~\ref{fig:ComputationTime_NrAbsentDrivers} shows the relation between the computational time and the number of absent drivers. We plot the median lines of the computational times for each instance. Overall, the computational time increases relatively slowly when we increase the number of absent drivers: when we double the number of absent drivers, the computational time approximately doubles. 

\begin{figure}[htb!]
\vspace{-4mm}
\centering
\includegraphics[scale=0.28]{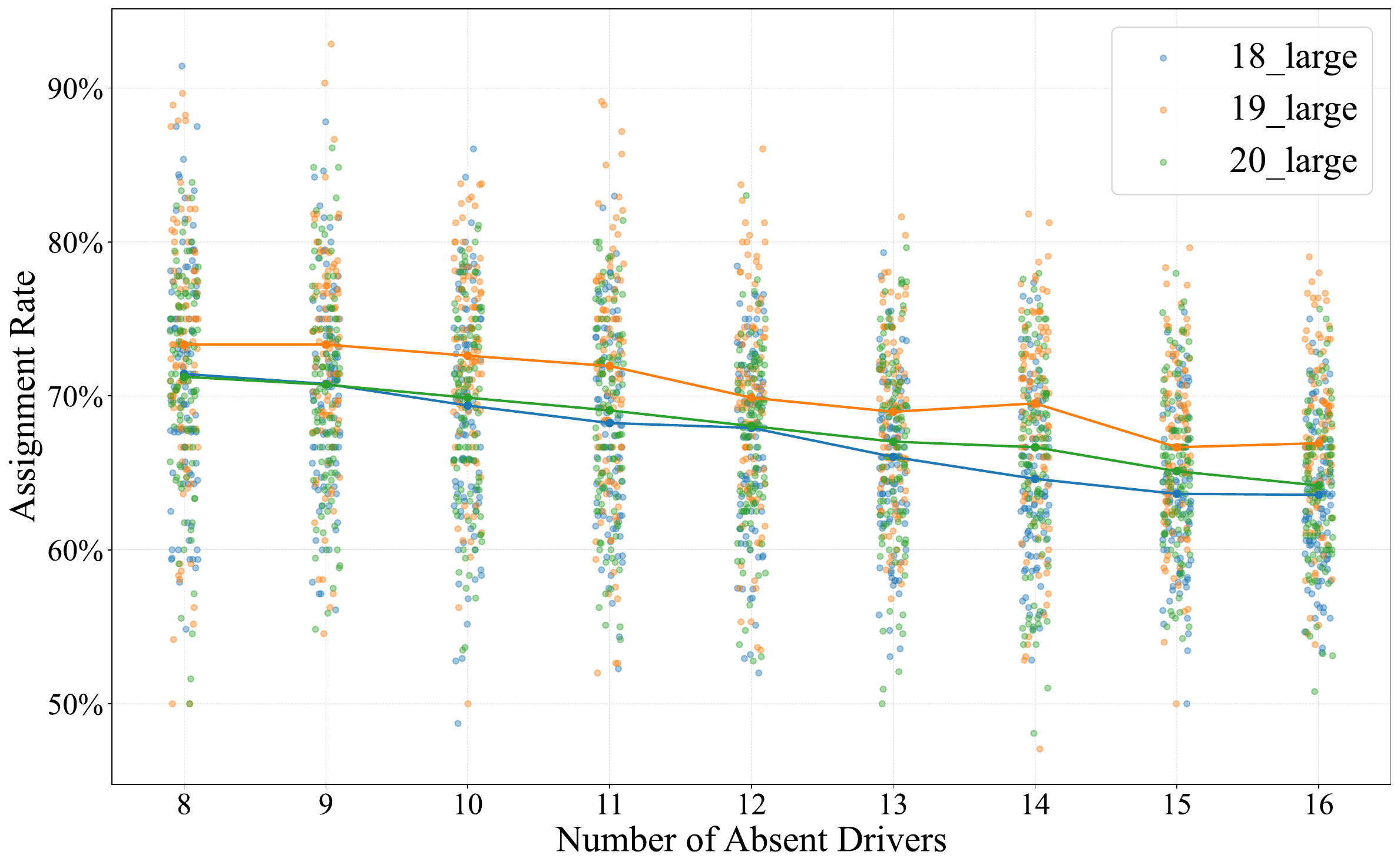}
\caption{Assignment rate distribution of the large instances.Points represent individual runs. Solid lines connect the median values computed for each number of absent drivers.}\label{fig:AssignedRate}
\vspace{-4mm}
\end{figure}

\begin{figure}[htb!]
\centering
\includegraphics[scale=0.28]{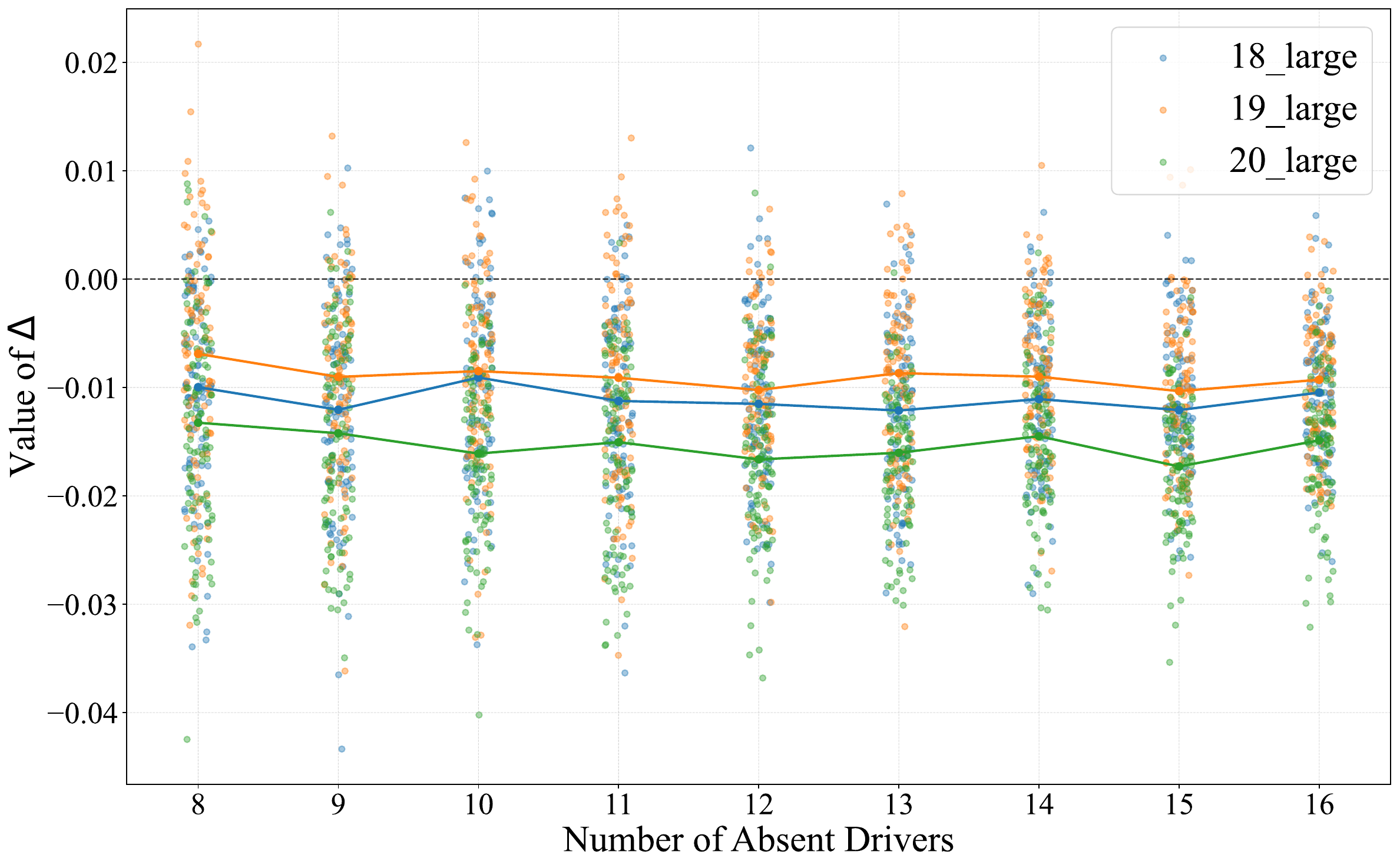}
\caption{$\Delta$ distribution of the large instances. Points represent individual runs. Solid lines connect the median values computed for each number of absent drivers. The dashed horizontal reference line marks $\Delta=0$.}
\label{fig:dwar}
\vspace{-4mm}
\end{figure}

Figure~\ref{fig:AssignedRate} displays the assignment rate distribution for the different large instances. The median value of the assignment rate varies between 60\% and 80\% when the number of absent drivers varies between 8 and 16. We obtain a higher assignment rate when the number of absent drivers decreases, which is reasonable as fewer tasks need to be assigned to more available drivers.

Figure~\ref{fig:dwar} reports the distribution of $\Delta$ values and the median lines for different numbers of absent drivers. The median values are all below zero, indicating that the remaining unassigned tasks are more difficult to assign than the assigned tasks. Although the point distribution varies across runs, the consistently negative median values show that this interpretation is stable under different numbers of absent drivers.

\subsubsection{Extension of Start and End Time Limit}\label{subsubsubsec:extension-tabu_time_adjust}

In this subsection, we examine the effect of relaxing drivers' start and end time limits. For each instance, we select one set of absent drivers whose original duties include either very early-starting tasks or very late-ending tasks, so that the rescheduling problem is sensitive to the imposed time limits. The extension varies from 0 to 3 hours, and we report the median results over 100 repeated runs for each setting in Table~\ref{tab:extension_start_end}. Larger extensions generally improve both the assignment rate and DWAR, indicating that additional temporal flexibility enables more unassigned tasks to be recovered. The main improvement occurs when the extension increases from 0 to 2 hours, whereas the additional gain from 2 to 3 hours is limited. The consistently negative values of $\Delta$ indicate that the remaining unassigned tasks are relatively difficult to assign. Computational times increase only moderately and remain low across all tested settings.

\begin{table}[htb!]
\centering
{\begin{tabular}{@{}cc|cccc@{}}
\hline
Instance  & \makecell[c]{extension\\hour} & \makecell[c]{assignment\\rate} & DWAR	& $\Delta$ & \makecell[c]{computational\\time}
\\ \hline
18\_large & 0 & 51.2\%& 48.7\%& -0.025& 3.8 s\\ \hline
18\_large & 1 & 56.1\% &53.5\% &-0.028&5.7 s\\ \hline
18\_large & 2 & 63.4\%	&60.5\%	&-0.028	&6.8 s\\ \hline
18\_large & 3 & 63.4\%	&60.6\%	&-0.029	&7.9 s\\ \hline
19\_large & 0 & 51.4\%	&49.2\%	&-0.022	&2.4 s\\ \hline
19\_large & 1 & 56.8\%	&54.5\%	&-0.024	&3.3 s\\ \hline
19\_large & 2 & 62.2\%	&60.3\%	&-0.024	&3.5 s\\ \hline
19\_large & 3 & 64.9\%	&62.5\%	&-0.024	&4.3 s\\ \hline
20\_large & 0 & 50.9\%	&59.6\%	&-0.013	&3.8 s\\ \hline
20\_large & 1 & 62.3\%	&60.2\%	&-0.021	&3.8 s\\ \hline
20\_large & 2 & 66.8\%	&64.9\%	&-0.019	&4.0 s\\ \hline
20\_large & 3 & 66.8\%	&65.0\%	&-0.019	&4.5 s\\ \hline
\end{tabular}}
\caption{Median results under start- and end-time extensions}
\label{tab:extension_start_end}
\end{table}

\section{Conclusions} \label{Conclusions}
In this paper, we considered a short-term crew rescheduling problem in which we aim to reschedule a one-day schedule the day before the operation because drivers' absences made the original plan infeasible. We formulated a tabu-search-based approach for our problem and adapted a column-generated approach as a benchmark. For both approaches, we use the same objective function and employ the same restrictions: operational limits (temporal and geographical consistency), bounds on the total working hours, restrictions from the drivers' licenses on vehicles and geographical areas, and bounds on rest times and on breaks. For the tabu-search-based approach, we built operationally feasible task connections in a task graph and used restricted graph search to generate feasible alternative paths.

The method was evaluated against exact and column-generation approach benchmarks. On reduced-scale instances, it produced solutions close to the optimum with substantially shorter computational times than the benchmarks. On full-scale instances, where solving the complete model directly is computationally intractable, it produced high-quality solutions within short runtimes. Supplementary ROMSOC-derived instances confirmed the computational advantage of our proposed method and highlighted the limitations of the diving approach of the column-generation benchmark. The graph-search evaluation showed that restricted DFS identifies feasible repair paths that cannot be represented by paths containing at most two task nodes. The analysis showed that the method remains computationally efficient as the number of absent drivers increases, and that additional operational flexibility through relaxed start and end time limits improves recovery performance.

For future work, we consider extending the proposed method in two main directions. Additional recovery actions, such as taxi-based deadheading and broader inter-operator deadheading, could be incorporated to improve schedule recovery.  Also, future extensions may include task priorities, driver preferences, and fairness considerations, allowing the rescheduling process to better reflect operational importance and human-centric requirements.

\section*{Acknowledgements} 
We would like to thank Olov Lindfeldt, Lars Swahn, Freddy Ahlström, and Viktor Reineck for the discussions and the input data from Mälartåg.


\clearpage
\appendix
\section{Adapted Column Generation}\label{CG}
In this section, we present our adapted formulation of the master problem and pricing problem from \citet{breugem2022column}. We use the column-generation approach as a benchmark for our tabu-search-based approach. In Subsection~\ref{RestrictedMasterProblem}, we introduce the adapted master and describe its difference from the original integer linear programming (ILP) model by \citet{breugem2022column}. In Subsection~\ref{PricingProblem}, we further display the adapted pricing problem.

\begin{figure}[htb!]
\centering
\includegraphics[scale=0.76]{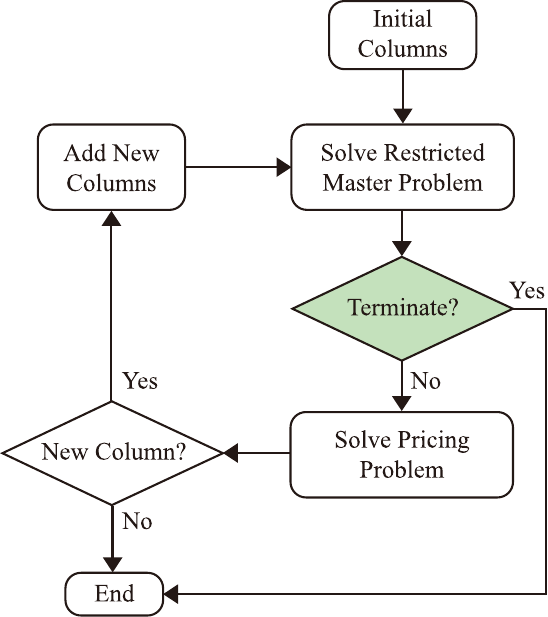}
\caption{The process of the adapted column-generation approach} 
    \label{fig:3}
\end{figure}

Figure~\ref{fig:3} shows the process of the column-generation approach. We start by solving the restricted master problem with the given initial columns, and then search for new columns by solving the pricing problem. We terminate the algorithm when no additional columns are generated. Additionally, during each assessment of the termination criteria, we verify whether the maximum allowed runtime has been exceeded. If this runtime limit is violated, we also terminate the algorithm.

Here, we present the model for the master problem and pricing problem that we adapted from \citet{breugem2022column} for short-term crew rescheduling. We adapt the ILP formulation of the column-generation approach and then simplify the master problem and pricing model. 

\subsection{Restricted Master Problem} \label{RestrictedMasterProblem}
\citet{breugem2022column} reschedule the tasks on both duty and roster levels. We only consider rescheduling tasks for a one-day schedule, so only duty-level rescheduling is needed. Therefore, we simplify the master problem by removing the variables and constraints related to rostering. We also change the pricing problem accordingly. Table~\ref{tab:1} defines the sets in the model. We also display the parameters and variables separately in Table~\ref{tab:2} and \ref{tab:3}, respectively.

\begin{table}[htb!]
\centering
\begin{tabular}{m{0.16\columnwidth}p{0.77\columnwidth}} \hline
\textbf{Sets} &\textbf{Description}\\ \hline
 $ G $              & Set of all tasks \\
 $ M $   &  Set of all available drivers including operating and standby drivers \\ 
 $K_m$ & Set of feasible duties for driver $m$, where $m \in M$\\
$K^{\mbox{\footnotesize{additional}}}$ & set of additional duties \\ \hline 
\end{tabular}
\caption{Sets}
    \label{tab:1}
\end{table}

\begin{table}[htb!]
\centering
\begin{tabular}{m{0.16\columnwidth}p{0.77\columnwidth}} \hline
 \textbf{Parameters}  & \textbf{Description} \\
 \hline
 $c_{k,m}$     &The cost of assigning duty $k$ for driver $m$\\
 $c_{k}^{\mbox{\footnotesize{additional}}}$     &The cost of additional duty $k$\\
 $a_{g,k}$     &Binary parameter, which equals 1 if task $g$ is included in duty $k$\\
\hline
\end{tabular}
\caption{Parameters}
    \label{tab:2}
\end{table}

\begin{table}[htb!]
\centering
\begin{tabular}{m{0.16\columnwidth}p{0.77\columnwidth}} \hline
 \textbf{Variables}  & \textbf{Description} \\
 \hline
 $ x_{k,m} $   & Binary variable, which equals 1 if duty $k$ is assigned to driver $m$, where $k \in K_m$, $m \in M$\\
$ y_{k} $   & Binary variable, which equals 1 if additional duty $k$ is selected, where $k \in K^{\mbox{\footnotesize{additional}}}$\\ \hline
\end{tabular}
\caption{Variables}
    \label{tab:3}
\end{table}

We present the restricted master problem. We aim to minimize the objective function~\eqref{obj}, which is a convex combination of the total cost of selected and additional duties.

\begin{flalign}
& \sum_{m \in M}\sum_{k \in K_m}{c_{k,m} x_{k,m}} + \sum_{k \in K^{\mbox{\footnotesize{additional}}}}{c_{k}^{\mbox{\footnotesize{additional}}}y_{k}}&  \label{obj}
\end{flalign}

Constraint~\eqref{con:coverall} ensures that each task $g$ is covered. Constraint~\eqref{con:assign} ensures that each driver $m$ is assigned one and only one duty. Constraints~\eqref{con:variable_x} and \eqref{con:variable_y} indicate the variable domains.
\begin{flalign}
& \sum_{m \in M}\sum_{k \in K_m}{a_{g,k} x_{k,m}} + \sum_{k \in K^{\mbox{\scriptsize{additional}}}}{a_{g,k}y_{k}} = 1&\forall g \in G  \label{con:coverall}
\end{flalign}
\begin{flalign}
& \sum_{k \in K_m}{x_{k,m}} = 1&\forall m \in M  \label{con:assign}
\end{flalign}
\begin{flalign}
& 0 \leq x_{k,m} \leq 1 &\forall k \in K_m,\quad m \in M  \label{con:variable_x}
\end{flalign}
\begin{flalign}
& 0 \leq y_{k} \leq 1 &\forall k \in K^{\mbox{\footnotesize{additional}}}  \label{con:variable_y}
\end{flalign}

\subsection{Pricing Problem} \label{PricingProblem}

Our reduced cost $RC_{k,m}$ is given by Equation~\eqref{con:ReduceCost}, where $\lambda_g$ is the dual variable of Constraint~\eqref{con:coverall}, and $\mu_m$ is the dual variable of Constraint~\eqref{con:assign}.
\begin{flalign}
& RC_{k,m} = c_{k,m} - \sum_{g \in G}{a_{g,k}\ \lambda_g} - \mu_m &  \label{con:ReduceCost}
\end{flalign}
\begin{flalign}
& RC^{\mbox{\footnotesize{additional}}}_{k} = c_{k}^{\mbox{\footnotesize{additional}}} - \sum_{g \in G}{a_{g,k}\ \lambda_g}&  \label{con:ReduceCost_1}
\end{flalign}

Both pricing problems are formulated as resource-constrained shortest path problems (RCSPPs) on duty graphs, and the constraint refers to the limitation of maximum duty length.

The proposed pricing adopts a two-phase strategy. A heuristic resource-constrained shortest path problem (RCSPP) is first applied to rapidly identify negative improving reduced-cost columns. The exact DFS is only activated when the heuristic phase fails to find any improving columns, and thereby prioritizes optimality over computing time.

\section{Difficulty-Weighted Assignment Rate}
\label{eva:difficulty_weighted_assignment_rate}
The standard assignment rate measures the proportion of recovered disrupted tasks but assumes all tasks are equally difficult to assign. This may not fully reflect the characteristics of the remaining unassigned tasks, since each task has different difficulty in being assigned. Therefore, we introduce the DWAR to better capture the assignment difficulties in this appendix.

Let $U^0$ denote the set of initially unassigned tasks and $U^f$ the set of tasks that remain unassigned after rescheduling. The standard assignment rate is defined as
$$
{\mathrm{AR}
=
1-\frac{|U^f|}{|U^0|}.}
$$

To account for heterogeneous task-assignment difficulty, each task $g \in U^0$ is assigned a positive weight $\omega_g$. The difficulty-weighted assignment rate is defined as
$$
{\mathrm{DWAR}
=
1-
\frac{
\sum_{g \in U^f} \omega_g
}{
\sum_{g \in U^0} \omega_g
}.}
$$

The weight is given by $$ \omega_g = 1+h_g, $$ where $h_g \in [0,1]$ is the normalized task-difficulty measurements. The measurement is computed as $$ h_g = \beta h_g^{G} + (1-\beta)h_g^{T}$$

$h_g^{G}$ denotes the graph-based difficulty component, which is calculated as the sum of the node degree of task $g$, including in-degree $e_g^{\mathrm{in}}$ and out-degree $e_g^{\mathrm{out}}$. The calculation function is displayed below.
$$
{h_g^{G}
=
\frac{1}{
\sqrt{
\left(e_g^{\mathrm{in}}+1\right)
\cdot
\left(e_g^{\mathrm{out}}+1\right)
}
}},
$$
A task with a low degree means that it has fewer valid arcs connected to it, which largely increases the difficulty of assigning this task to the schedule. Also, the temporal component $h_g^{T}$ is defined as follows:

$$
{h_g^{T}
=
1-2\cdot \min\{\tau_g,1-\tau_g\} \text{, where }\tau_g =
\frac{t_g - T_{\min}}{T_{\max}-T_{\min}}}$$

The symbol $t_g$ denotes the start time of task $g$. The term $\min\{\tau_g,1-\tau_g\}$ measures the normalized distance from task $g$ to the closer boundary of the whole planning horizon, which starts from $T_{\min}$ to $T_{\max}$. Hence, $h_g^{T}$ increases as task $g$ moves closer to the beginning or end of the horizon: it is equal to 0 at the middle of the horizon and reaches 1 at either boundary. A task that is closer to the horizon boundaries compared to other tasks typically has fewer available drivers to select, which also largely increases the assignment difficulty. DWAR thus complements the assignment rate by indicating whether the remaining unassigned tasks are relatively easy or difficult to recover. In addition, we set $\beta=0.7$ in our experiments.

Finally, we calculate the difference $\Delta$ between the difficulty-weighted assignment rate $\mathrm{DWAR}$ and the assignment rate $\mathrm{AR}$.

$$
{\Delta
=
\mathrm{DWAR}-\mathrm{AR}.}
$$
The difference, $\Delta$, provides a simple metric, measuring the assignment difficulty of the remaining unassigned tasks. A negative $\Delta$ indicates that the remaining unassigned tasks are relatively more difficult to assign. DWAR and $\Delta$ are used only for analysis and are not included in the optimization objective.

Here, we give an example of computing the value of DWAR and $\Delta$. For example, suppose that two tasks are initially unassigned:
$$
{U^0=\{1,2\}.}
$$
Their difficulty weights are
$$
{\omega_1=1.175,\qquad \omega_2=1.364.}
$$
After rescheduling, task 1 is recovered, while task 2 remains unassigned. Then, we have $U^f=\{2\}$.

$$
{\mathrm{DWAR}
=
1-
\frac{
\sum_{g \in U^f} \omega_g
}{
\sum_{g \in U^0} \omega_g
}
=
1-\frac{1.364}{1.175+1.364}
= 0.463
}
$$

$$
{\Delta
=
DWAR-AR
=0.463-0.5
= -0.037
}
$$

The negative value of $\Delta$ indicates that the remaining unassigned task has a higher difficulty weight than the already recovered task.

\end{document}